\documentclass[letterpaper, 10 pt, conference]{ieeeconf}

\usepackage{amsmath}
\usepackage{graphicx}
\usepackage{amssymb}
\usepackage{multirow}
\usepackage{color}

\newtheorem{theorem}{Theorem}

\newtheorem{corollary}{Corollary}

\title{The Barrier Surface in the Cooperative Football Differential Game}
\author{Eloy Garcia, David W. Casbeer, and Meir Pachter
\thanks{E. Garcia and D. Casbeer are with the Control Science Center of Excellence, Air Force Research Laboratory, Wright-Patterson AFB, OH 45433. Corresponding author \ttfamily{eloy.garcia.2@us.af.mil}}
\thanks{M. Pachter is with the Department of Electrical Engineering, Air Force Institute of Technology, Wright-Patterson AFB, OH 45433.}
}
          
\begin{document}
\maketitle 

\begin{abstract}
This paper considers the blocking or football pursuit-evasion differential game. Two pursuers cooperate and try to capture the ball carrying evader as far as possible from the goal line. The evader wishes to be as close as possible to the goal line at the time of capture and, if possible, reach the line. In this paper the solution of the game of kind is provided: The Barrier surface that partitions the state space into two winning sets, one for the pursuer team and one for the evader, is constructed.  Under optimal play, the winning team is determined by evaluating the associated Barrier function.
\end{abstract}

     

\section{Introduction} \label{sec:intro}
Determining the outcome of the conflict ahead of time is a one of the most important problems in pursuit-evasion differential games. This is commonly known as solving the game of kind and it provides the answer to which team or player wins the game under optimal play given the current state.
Pursuit-evasion conflicts are representative of many important but challenging problems in aerospace, control, and robotics \cite{Weintraub20}. Many different pursuit-evasion conflicts were formulated in the seminal work by Isaacs \cite{Isaacs65}. Two important concepts were introduced in \cite{Isaacs65}: the game of kind and the game of degree. The solution to the latter provides the actual players' strategies that achieve the prescribed outcome provided by the solution of the game of kind.



Recent authors have investigated different approaches to pursuit-evasion games which rely on geometric constructs    \cite{zhou2016cooperative,pierson2016intercepting} or using discrete-time analysis and matrix games \cite{Selvakumar17,Casini19}. A differential game framework is generally desired \cite{chen2016multiplayer}, but is often avoided due to the perceived challenges in solving the Hamilton-Jacobi-Isaacs (HJI) equation \cite{pierson2016intercepting}.
In reach and avoid games \cite{margellos2011hamilton}, the HJI is solved numerically and computation time significantly increases as the state dimension increases.
Reach and avoid games where an evader strives to reach a goal set while the pursuer tries to intercept the evader before it reaches the goal, include the capture the flag game  \cite{Blake04,Huang15,GarciaCDC18}, defending a moving target \cite{fisac2015pursuit,Coon17,Garcia2019,harini2015new,Weintraub2018} or static asset \cite{Shishika18,Garcia20LCS,YanShi17}, and assisting and rescuing teammates  \cite{Oyler16,Scott13}.

In this paper we revisit the football differential game where two cooperating pursuers try to capture an evader who aims to reach the goal line. 
This problem was addressed in \cite{GarciaACC19}, but only the game of degree was solved; it was assumed that the game starts in the region of win of the pursuers. This problem has immediate applications to border defense using autonomous vehicles \cite{Bopardikar11,Vonmoll19}. 
In the present paper the game of kind is solved; the solution provides the answer to which team wins the game under optimal play. Additionally, when the pursuer team wins  we  provide the mode of capture in the optimal outcome: single pursuer capture or dual  pursuer cooperative capture. 
The game of kind has been recently addressed in \cite{YanShi19} which only considered the case where the pursuers are faster than the evader but the pursuers have the same speed between them. This paper offers several improvements compared to \cite{YanShi19}:
 In this paper we consider two scenarios: 1) the case where the evader is not slower than the pursuers but it has the same speed and 2) the case where the pursuers are faster. In the latter, we generalize to the case where the pursuers have different speeds.
We provide explicit expressions of the Barrier surface and also of the cross-sections of the Barrier surface. We do this by defining the type of curves while \cite{YanShi19} used less explicit `base curves'; these `base curves' only describe conditions for capture at the goal line without explicitly describing their geometric properties. Similarly, we determine the Barrier surface which is a surface on a higher dimensional state space. For instance, for two pursuers and one evader, the dimension of the state space is six. Then, a cross-section of the Barrier surface can be projected into the 2D plane by fixing the positions of the pursuers. 
The main advantage with respect to \cite{YanShi19} is that we solve for an explicit expression of the Barrier surface without the use of `base curves'. This feature allows for simpler and more explicit expressions that highlight their geometric properties.

The paper is organized as follows. The problem is formulated in Section \ref{sec:Problem}. The Barrier surface when all players have the same speed is obtained in Section \ref{sec:ss} and for fast pursuers with different speeds in Section \ref{sec:FP}. 
Conclusions are drawn in Section \ref{sec:concl}.


\section{The Football Pursuit-Evasion DG} \label{sec:Problem}

Consider the three agents $E$, $P_1$, and $P_2$ with max speeds given by $v_E$, $v_1$, and $v_2$. The states of $E$, $P_1$, and $P_2$ are $\textbf{x}_E=(x_E,y_E)$, $\textbf{x}_1=(x_1,y_1)$, and $\textbf{x}_2=(x_2,y_2)$. Let 
\begin{align}
\textbf{x}:=( x_E, y_E, x_1, y_1, x_1, y_1)\in \mathbb{R}^6
\end{align}
 denote the complete state of the differential game.
The control variable of $E$ is $\textbf{u}_A=\left\{\phi\right\}$, its heading angle. Players $P_1$ and $P_2$ form a team and they cooperatively choose their respective headings $\textbf{u}_B=\left\{\psi_1,\psi_2\right\}$. 
The dynamics $\dot{\textbf{x}}=\textbf{f}(\textbf{x},\textbf{u}_A,\textbf{u}_B)$ are specified by the system of differential equations
\begin{equation}
\begin{alignedat}{2}
	\dot{x}_E&=v_E\cos\phi,  &\qquad x_E(0)&=x_{E_0}      \\
	\dot{y}_E&=v_E\sin\phi,   &\qquad  y_E(0)&=y_{E_0}    \\
        \dot{x}_1&=v_1\cos\psi_1,   &\qquad x_1(0)&=x_{1_0}   \\
	\dot{y}_1&=v_1\sin\psi_1,  &\qquad y_1(0)&=y_{1_0}  \\
	\dot{x}_2&=v_2\cos\psi_2,   &\qquad x_2(0)&=x_{2_0}  \\
	\dot{y}_2&=v_2\sin\psi_2,   &\qquad y_2(0)&=y_{2_0}   \label{eq:xT}
\end{alignedat}
\end{equation}
where the admissible controls are the players' headings $\phi,\psi_1,\psi_2 \in [-\pi,\pi]$ and 
\begin{align}
\textbf{x}_0 := (x_{E_0}, y_{E_0}, x_{1_0}, y_{1_0}, x_{2_0}, y_{2_0}) = \textbf{x}(t_0)
\end{align}
 is the initial state.
The game is played in a rectangular region $\Omega \subset\mathbb{R}^2$ of the Cartesian plane and without loss of generality we have $\Omega:= \big\{ x,y | x\in[0,\bar{x}], y\geq0 \big\} $, where $\bar{x}$ is a specified parameter. The goal line is  $G:= \big\{ x,y\in \Omega | y=0 \big\}$. 
With respect to the information constraints, every agent knows the dynamics \eqref{eq:xT}. It is assumed that the agents use causal strategies and that every agent has access to the state $\textbf{x}$ at the current time $t$, that is, the capture game is a perfect information differential game; the optimal strategies will be state feedback strategies. Finally, and most importantly, we assume that the agents do not know the opponent's current decision, we assume non-anticipative strategies.



The game terminates if $y_E=0$ and the evader wins by reaching the goal line $G$ before being captured by any of the pursuers. The game also terminates if any of the pursuers captures the evader, where point capture is considered. The terminal set is then given by
\begin{align}
   \left.
	 \begin{array}{l l}
	  \mathcal{T}:=     \big\{ \ \textbf{x} \ | y_E =0 \big\} \cup  \big\{ \ \textbf{x} \ | x_1=x_E, y_1=y_E \big\} \\
   \qquad \quad \cup  \big\{ \ \textbf{x} \ | x_2=x_E, y_2=y_E \big\}.   \label{eq:Set}
\end{array}  \right. 
\end{align}
%
We refer to the pursuer capturing the evader as the active pursuer. If simultaneous capture occurs then both pursuers are active. Thus, a complete solution of the game of kind will be achieved by determining capture conditions and the pursuer which will capture the evader under optimal play.

In reference \cite{GarciaACC19} the game of degree was addressed assuming the state of the game is in the region of win of the pursuers, that is, the pursuers are guaranteed to capture the evader, provided they play optimlly. The players try to minmax the terminal separation from the capture point to the goal line. However, the Barrier surface was not constructed and it was not possible in \cite{GarciaACC19} to evaluate if capture is in fact possible. This paper fully addresses this problem by obtaining the analytical Barrier surface that separates the parties regions of win in the state space.
 The state space $\mathbb{R}^6$ is partitioned into two sets: $\mathcal{R}_p$ and $\mathcal{R}_e$ which are defined as follows
\begin{align}
\left.
	 \begin{array}{l l}
\mathcal{R}_p:=  \big\{ \ \textbf{x} \ |  \ B( \textbf{x})<0  \big\} \\
\mathcal{R}_e:=  \big\{ \ \textbf{x} \ |  \  B( \textbf{x})>0  \big\}
\end{array}   \right.  \label{eq:ReDef}
\end{align}
and the Barrier surface, which separates the two sets $\mathcal{R}_p$ and $\mathcal{R}_e$, is specified by
\begin{align}
\left.
	 \begin{array}{l l}
\mathcal{B}:=  \big\{ \ \textbf{x} \ | \ B( \textbf{x})=0  \big\} \\
\end{array}   \right.  \label {eq:BarrSurface}
\end{align}
where the Barrier function $B( \textbf{x})$ is explicitly obtained for each of the two cases analyzed in this paper in Sections \ref{sec:ss} and \ref{sec:FP}.

\section{Same Speed Players}  \label{sec:ss}
\begin{figure}
	\begin{center}
		\includegraphics[width=8.4cm,trim=4.0cm 1.9cm 1.7cm .5cm]{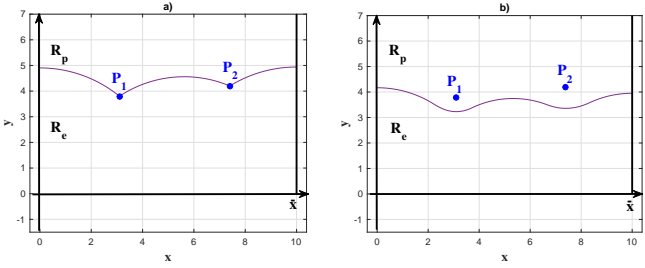}
	\caption{Examples of cross-section of Barrier surface: a) Same speed players. b) Fast pursuers}
	\label{fig:ex}
	\end{center}
\end{figure}

The case where all agents have the same speed is considered first.
We assume, without loss of generality, that $x_1<x_2$, that is, $P_2$ is located to the right of $P_1$ in Fig. \ref{fig:ex}. 


\begin{theorem}   \label{th:sstheo}
Consider the cooperative football differential game with two pursuers and one evader. The Barrier surface is given by $B(\textbf{x})=0$ where $B(\textbf{x})$ consists of three different segments
\begin{align}
	B(\textbf{x})= \left\{ 
	 \begin{array}{l l}
	B_1(\textbf{x}), \ \ \  \text{if} \ x_E\leq x_1 \\
         B_2(\textbf{x}), \ \ \ \text{if} \  x_1<x_E < x_2 \\
	 B_3(\textbf{x}), \ \ \  \text{if} \  x_E\geq x_2
\end{array}  \right.  \label{eq:Barriersurf}
\end{align}
each segment is explicitly given in terms of the state by
\begin{align}
 \left.
	 \begin{array}{l l}
	B_1(\textbf{x})= x_1^2+y_1^2-x_E^2-y_E^2  \\
	B_2(\textbf{x})=  (x_1-x_I)^2+y_1^2-(x_E-x_I)^2-y_E^2 \\
	 B_3(\textbf{x})= (x_2-\bar{x})^2+y_2^2-(x_E-\bar{x})^2-y_E^2
\end{array}  \right.  \label{eq:Barrier}
\end{align}
where $x_I=\frac{1}{2} \frac{x_2^2+y_2^2-x_1^2-y_1^2}{x_2-x_1}$. 
\end{theorem}
\textit{Proof}. The dominance regions between two players with the same speed are separated by the orthogonal bisector (OBS) of the line connecting the positions of the players. Consider the first case where $x_E\leq x_1$. In this case, only $P_1$ is active and the OBS of the segment $\overline{P_1E}$ is given by the line $y=m_1x+n_1$ where
 \begin{align}
 \left.
	 \begin{array}{l l}
 m_1=-\frac{x_1-x_E}{y_1-y_E} \\
 n_1=\frac{1}{2}\frac{x_1^2+y_1^2-x_E^2-y_E^2}{y_1-y_E}.
 \end{array}  \right.   \label{eq:m1n1}
 \end{align}
  The Barrier surface is obtained when  the line $y=m_1x+n_1$ intersects the $x$-axis at $x=0$. Making these substitutions we have
\begin{align}
 \left.
	 \begin{array}{l l}
	&0=m_1(0) +n_1 \\
	 \Rightarrow &n_1=0.
\end{array}  \right.   \label{eq:OBS1}
\end{align}
Substituting \eqref{eq:m1n1} into \eqref{eq:OBS1} we obtain that the Barrier surface is $B_1(\textbf{x})=0$ where $B_1(\textbf{x})$ is given in \eqref{eq:Barrier}.

The second case is when $x_1<x_E < x_2$ and both pursuers are active. The OBS of the segment $\overline{P_1P_2}$ is given by the line $y=mx+n$, where 
 \begin{align}
 \left.
	 \begin{array}{l l}
m=-\frac{x_2-x_1}{y_2-y_1}  \\
n=\frac{1}{2}\frac{x_2^2+y_2^2-x_1^2-y_1^2}{y_2-y_1}.
\end{array}  \right.   \label{eq:mn}
\end{align}
Let the intersection point between the line $y=mx+n$ and the $x$-axis be denoted by $I=(x_I,0)$ where $x_I$ is obtained as follows
\begin{align}
 \left.
	 \begin{array}{l l}
	&0=mx_I +n \\
	\Rightarrow &x_I=-\frac{n}{m} = \frac{1}{2} \frac{x_2^2+y_2^2-x_1^2-y_1^2}{x_2-x_1}.
\end{array}  \right.  \label{eq:OBSp}
\end{align}
Simultaneous capture of $E$ by $P_1$ and $P_2$ occurs at the intersection of the OBS of the segment $\overline{P_1E}$ and the OBS of the segment $\overline{P_2E}$, where the latter is given by $y=m_2x+n_2$ where 
 \begin{align}
 \left.
	 \begin{array}{l l}
m_2=-\frac{x_2-x_E}{y_2-y_E} \\
n_2=\frac{1}{2}\frac{x_2^2+y_2^2-x_E^2-y_E^2}{y_2-y_E}.
\end{array}  \right.   \label{eq:m2n2}
\end{align}
Since the players have the same speed all lines $y=mx+n$,  $y=m_1x+n_1$, and $y=m_2x+n_2$  intersect at the same point.
The Barrier surface is obtained when the intersection point is given by $I=(x_I,0)$; then, we have that
\begin{align}
 \left.
	 \begin{array}{l l}
	0=m_1x_I +n_1  \\
	\ \ =- 2 (x_1-x_E) x_I + x_1^2+y_1^2-x_E^2-y_E^2
\end{array}  \right.  \label{eq:sscase2}
\end{align}
adding and subtracting $x_I^2$ to \eqref{eq:sscase2} we have the Barrier surface $B_2(\textbf{x})=0$ where $B_2(\textbf{x})$ is given in \eqref{eq:Barrier}.

Finally, the third case is when $x_E\geq x_2$ and only $P_2$ is active. The Barrier surface is obtained when the line $y=m_2x+n_2$ intersects the point $(\bar{x},0)$
\begin{align}
 \left.
	 \begin{array}{l l}
	0=m_2 \bar{x} +n_2  \\
	\ \ =- 2 (x_2-x_E) \bar{x} + x_2^2+y_2^2-x_E^2-y_E^2
\end{array}  \right.  \nonumber  
\end{align}
and the Barrier surface is $B_3(\textbf{x})=0$ where $B_3(\textbf{x})$ is given in \eqref{eq:Barrier}. \ $\square$

Note that the Barrier surface $B(\textbf{x})=0$  is a surface in the sixth-dimensional state space $\mathbb{R}^6$. In the first and third cases only one of the pursuers is active and the state of the other pursuer is not relevant. In the second case all states are necessary, note that $x_I$ is a function of the states of both pursuers.
It is also possible to address a particular case, as it was done in \cite{YanShi19}, where the positions of the pursuers are fixed; then, a \textit{cross-section} of the Barrier surface is obtained and can be illustrated in the 2D plane. 

\begin{corollary}   \label{cor:sscoro}
For given positions of $P_1$ and $P_2$, the cross-section of the Barrier surface is given by arc segments corresponding to three different circles
 \begin{align}
 \left.
	 \begin{array}{l l}
	x^2+y^2 = x_1^2+y_1^2, \ \ \ \ \  \text{if} \ x\leq x_1    \\
	(x\!-\!x_I)^2+y^2 =  (x_1\!-\!x_I)^2+y_1^2,  \ \text{if} \  x_1<x < x_2   \\
	(x\!-\!\bar{x})^2+y^2  =(x_2\!-\!\bar{x})^2+y_2^2, \ \ \  \text{if} \  x\geq x_2
\end{array}  \right.  \nonumber  
\end{align}
\ $\square$
\end{corollary}
Figure \ref{fig:ex}.a shows an example of the cross-section of the Barrier surface for given positions of the pursuers. The projection of the sets $\mathcal{R}_e$ and $\mathcal{R}_p$ into the 2D plane is also illustrated in Fig. \ref{fig:ex}.a; they are denoted, respectively, by $R_e$ and $R_p$. If the evader is initially located in $R_p$ then, under optimal play, it will be captured by the pursuer team. However, if the evader is initially located in $R_e$ then it will be able to reach the goal line before being captured by any of the pursuers.

\section{Faster Pursuers}  \label{sec:FP}
In this section the case where the evader is slower than the pursuers is addressed.
Define the speed ratios $\gamma_1=v_E/v_1<1$ and $\gamma_2=v_E/v_2<1$. We consider in general the case where  $\gamma_1\neq \gamma_2$, that is, the pursuers have different speeds. We assume, without loss of generality, that $P_2$ is faster than $P_1$, hence $\gamma=v_1/v_2<1$.


\begin{theorem}   \label{th:FPtheo}
Consider the cooperative football differential game with two fast pursuers and one evader. The Barrier surface is given by $B(\textbf{x})=0$ where $B(\textbf{x})$ consists of five different segments
\begin{align}
 \left.
	 \begin{array}{l l}
	B(\textbf{x})= \\
	 \left\{ 
	 \begin{array}{l l}
	B_1(\textbf{x}), \  \text{if} \ x_E \leq \gamma_1^2 x_1 \\
         B_2(\textbf{x}), \  \text{if} \  \gamma_1^2 x_1<x_E \leq (1\!-\!\gamma_1^2)x_I + \gamma_1^2 x_1 \\
	 B_3(\textbf{x}), \  \text{if} \  (1\!-\!\gamma_1^2)x_I \!+\! \gamma_1^2 x_1< x_E  \leq(1\!-\!\gamma_2^2)x_I \!+\! \gamma_2^2 x_2  \\
	 B_4(\textbf{x}), \  \text{if} \  (1\!-\!\gamma_2^2)x_I \!+\! \gamma_2^2 x_2< x_E  \leq (1\!-\!\gamma_2^2)\bar{x} \!+\! \gamma_2^2 x_2  \\
	 B_5(\textbf{x}), \  \text{if} \  x_E  > (1\!-\!\gamma_2^2)\bar{x} + \gamma_2^2 x_2      
\end{array}  \right.  
\end{array}  \right.  \label{eq:BarriersurfFP}
\end{align}
each segment is explicitly given in terms of the state by
\begin{align}
 \left.
	 \begin{array}{l l}
	B_1(\textbf{x})=\gamma_1^2(x_1^2+y_1^2)-x_E^2-y_E^2  \\
	B_2(\textbf{x})=  \gamma_1^2\big(x_1^2+(1\!-\!\gamma_1^2)y_1^2\big) +\gamma_1^2x_E^2 - (1\!-\!\gamma_1^2)y_E^2   \\ 
	 \qquad  \qquad  -2\gamma_1^2x_1x_E   \\
	 B_3(\textbf{x})=   \gamma_1^2\big((x_1-x_I)^2+y_1^2\big)- (x_E-x_I)^2-y_E^2  \\
	 B_4(\textbf{x})=  \gamma_2^2\big(x_2^2+(1\!-\!\gamma_2^2)y_2^2\big) +\gamma_2^2x_E^2 - (1\!-\!\gamma_2^2)y_E^2  \\
	 \qquad  \qquad  -2\gamma_2^2x_2x_E   \\
	 B_5(\textbf{x})=\gamma_2^2\big((x_2-\bar{x})^2+y_2^2\big)-(x_E-\bar{x})^2-y_E^2
\end{array}  \right.  \label{eq:BarrierFP}
\end{align}
where now 
\begin{align}
 \left.
	 \begin{array}{l l}
x_I=\frac{ x_1-\gamma^2x_2 + \sqrt{\gamma^2(x_1-x_2)^2 - (1-\gamma^2)(y_1^2 - \gamma^2y_2^2) } }{1-\gamma^2}. 
\end{array}  \right.  \label{eq:XIFP}
\end{align}
\end{theorem}

\textit{Proof}. Since the pursuers are faster than the evader, now the dominance regions between each pair of players are separated by an Apollonius circle. The Apollonius circle between $E$ and $P_1$ is given by
\begin{align}
 \left.
	 \begin{array}{l l}
	(x-a_1)^2 + (y-b_1)^2 = r_1^2
\end{array}  \right.  \label{eq:ACP1E}
\end{align}
where 
\begin{align}
 \left.
	 \begin{array}{l l}
a_1=\frac{1}{1-\gamma_1^2}(x_E-\gamma_1^2x_1) \\
b_1=\frac{1}{1-\gamma_1^2}(y_E-\gamma_1^2y_1) \\
r_1=\frac{\gamma_1}{1-\gamma_1^2}\sqrt{(x_E\!-\!x_1)^2+(y_E\!-\!y_1)^2}.
\end{array}  \right.  \nonumber 
\end{align}
 When $x_E \leq \gamma_1^2 x_1$ only $P_1$ is active and the Barrier surface is obtained when the $EP_1$ Apollonius circle intersects the $x$-axis at $x=0$. Hence,  \eqref{eq:ACP1E} becomes $a_1^2 + b_1^2 = r_1^2$. We can write this expression explicitly in terms of the state as follows
\begin{align}
 \left.
	 \begin{array}{l l}
	 \gamma_1^2(1-\gamma_1^2)(x_1^2+y_1^2)-(1-\gamma_1^2)(x_E^2+y_E^2)  = 0
\end{array}  \right.  \nonumber
\end{align}
cancelling the common factor $(1-\gamma_1^2)$ we obtain the Barrier surface $B_1(\textbf{x})=0$ where $B_1(\textbf{x})$ is given in \eqref{eq:BarrierFP}.

The second case is when $\gamma_1^2 x_1<x_E \leq (1\!-\!\gamma_1^2)x_I + \gamma_1^2 x_1$ holds. In this case $P_1$ is still the only active pursuer but now the Barrier surface is obtained when the $EP_1$ Apollonius circle is tangent to the $x$-axis. The point of tangency is such that $x=a_1$ and $y=0$. From \eqref{eq:ACP1E}, we then have that $b_1-r_1=0$; this is equivalent to 
\begin{align}
 \left.
	 \begin{array}{l l}
	\quad \ \frac{y_E - \gamma_1^2y_1 -\gamma_1\sqrt{(x_E\!-\!x_1)^2+(y_E\!-\!y_1)^2}}{1-\gamma_1^2}   = 0  \\
	 \Rightarrow y_E^2 + \gamma_1^4y_1^2 = \gamma_1^2 (x_E^2+y_E^2 +x_1^2+y_1^2 -2x_1x_E).
\end{array}  \right.  \nonumber  
\end{align}
Rearranging the terms of the previous expression we obtain the Barrier surface $B_2(\textbf{x})=0$ where $B_2(\textbf{x})$ is given in \eqref{eq:BarrierFP}.

Simultaneous capture of $E$ by both pursuers occurs when  $(1\!-\!\gamma_1^2)x_I \!+\! \gamma_1^2 x_1< x_E  \leq(1\!-\!\gamma_2^2)x_I \!+\! \gamma_2^2 x_2$ and $E$ is captured at the intersection of the $EP_1$ and the $EP_2$ Apollonius circles. The latter is given by 
\begin{align}
 \left.
	 \begin{array}{l l}
	(x-a_2)^2 + (y-b_2)^2 = r_2^2
\end{array}  \right.  \label{eq:ACP2E}
\end{align}
where 
\begin{align}
 \left.
	 \begin{array}{l l}
a_2=\frac{1}{1-\gamma_2^2}(x_E-\gamma_2^2x_2) \\
b_2=\frac{1}{1-\gamma_2^2}(y_E-\gamma_2^2y_2)  \\
r_2=\frac{\gamma_2}{1-\gamma_2^2}\sqrt{(x_E\!-\!x_2)^2+(y_E\!-\!y_2)^2}.
 \end{array}  \right.  \nonumber   
\end{align}
We also note that a third circle, the $P_1P_2$ Apollonius circle, intersects the other two at the same intersection points. The $P_1P_2$ Apollonius circle is given by
\begin{align}
 \left.
	 \begin{array}{l l}
	(x-a)^2 + (y-b)^2 = r^2
\end{array}  \right.  \label{eq:ACPP}
\end{align}
where 
\begin{align}
 \left.
	 \begin{array}{l l}
a=\frac{1}{1-\gamma^2}(x_1-\gamma^2x_2) \\
b=\frac{1}{1-\gamma^2}(y_1-\gamma^2y_2)   \\
r=\frac{\gamma}{1-\gamma^2}\sqrt{(x_1\!-\!x_2)^2+(y_1\!-\!y_2)^2}.
\end{array}  \right.  \nonumber
\end{align} 
The Barrier surface is obtained when the intersection point is given by $I=(x_I,0)$; then, using \eqref{eq:ACPP} we have
\begin{align}
 \left.
	 \begin{array}{l l}
	x^2-2ax+a^2+b^2-r^2=0
\end{array}  \right.  \label{eq:ACinter}
\end{align}
where the applicable solution of \eqref{eq:ACinter} is 
\begin{align}
 \left.
	 \begin{array}{l l}
	x_I = a + \sqrt{r^2-b^2}
\end{array}  \right.  \label{eq:xIFP-2}
\end{align}
substituting the values of $a$, $b$, and $r$ into \eqref{eq:xIFP-2} we obtain $x_I$  explicitly in terms of the state as given in \eqref{eq:XIFP}. 
To determine the Barrier surface we can use \eqref{eq:XIFP} and the  $EP_1$ Apollonius circle so we have $(x_I-a_1)^2+b_1^2=r_1^2$ which can be written as follows
\begin{align}
 \left.
	 \begin{array}{l l}
	 \big( (1-\gamma_1)x_I - (x_E-\gamma_1^2x_1) \big)^2  + (y_E-\gamma_1^2y_1)^2  \\
	  = \gamma_1^2  \big(  (x_E-x_1)^2 + (y_E-y_1)^2 \big)^2
\end{array}  \right.  \nonumber
\end{align}
which can be written as
\begin{align}
 \left.
	 \begin{array}{l l}
       \gamma_1^2 (x_1^2+y_1^2) + 2x_I(x_E-\gamma_1^2x_1) \\
       -(1-\gamma_1^2)x_I^2   -x_E^2 -y_E^2 =0  
\end{array}  \right.  \nonumber
\end{align}
and, rearranging terms, we obtain the Barrier surface $B_3(\textbf{x})=0$ where $B_3(\textbf{x})$ is given in \eqref{eq:BarrierFP}.

The fourth case occurs when $(1\!-\!\gamma_2^2)x_I \!+\! \gamma_2^2 x_2< x_E  \leq (1\!-\!\gamma_2^2)\bar{x} \!+\! \gamma_2^2 x_2$ and only $P_2$ is active. The Barrier surface is obtained when the $EP_2$ Apollonius circle is tangent to the $x$-axis. This case is similar to the second case and the Barrier surface is $B_4(\textbf{x})=0$ where $B_4(\textbf{x})$ is given in \eqref{eq:BarrierFP}. 
Finally, when $x_E  > (1\!-\!\gamma_2^2)\bar{x} + \gamma_2^2 x_2$, the Barrier surface is obtained when the $EP_2$ Apollonius circle intersects the $x$-axis at $x=\bar{x}$. This case is similar to the first case and the Barrier surface is $B_5(\textbf{x})=0$ where $B_5(\textbf{x})$ is given in \eqref{eq:BarrierFP}. \ $\square$

Similar to the case of same speed players, the Barrier surface $B(\textbf{x})=0$  is a surface in the sixth-dimensional state space $\mathbb{R}^6$ and we can fix the positions of the pursuers in order to obtain a \textit{cross-section} of the Barrier surface which can be illustrated in the 2D plane. 

\begin{corollary}   \label{cor:FPcoro}
For given positions of $P_1$ and $P_2$, the cross-section of the Barrier surface is given by arc segments corresponding to three different circles and the arc segments of two different hyperbolas
 \begin{align}
 \left.
	 \begin{array}{l l}
	x^2+y^2 =  \gamma_1^2(x_1^2+y_1^2)   \\
	   \qquad  \qquad    \text{if} \ x \leq \gamma_1^2 x_1;  \\
         \gamma_1^2x^2 - (1\!-\!\gamma_1^2)y^2   -2\gamma_1^2x_1x + \gamma_1^2\big(x_1^2+(1\!-\!\gamma_1^2)y_1^2\big)  =0   \\
         \qquad  \qquad   \text{if} \  \gamma_1^2 x_1<x \leq (1\!-\!\gamma_1^2)x_I + \gamma_1^2 x_1; \\
	(x-x_I)^2 + y^2 = \gamma_1^2\big((x_1-x_I)^2+y_1^2\big) \\  
	\qquad \qquad    \text{if} \  (1\!-\!\gamma_1^2)x_I \!+\! \gamma_1^2 x_1< x  \leq(1\!-\!\gamma_2^2)x_I \!+\! \gamma_2^2 x_2;  \\
	 \gamma_2^2x^2 - (1\!-\!\gamma_2^2)y^2   -2\gamma_2^2x_2x + \gamma_2^2\big(x_2^2+(1\!-\!\gamma_2^2)y_2^2\big)  =0   \\
         \qquad  \qquad  \text{if} \  (1\!-\!\gamma_2^2)x_I \!+\! \gamma_2^2 x_2< x  \leq (1\!-\!\gamma_2^2)\bar{x} \!+\! \gamma_2^2 x_2;  \\
	 (x-\bar{x})^2+y^2 =\gamma_2^2\big((x_2-\bar{x})^2+y_2^2\big) \\
	\qquad \qquad    \text{if} \  x> (1\!-\!\gamma_2^2)\bar{x} + \gamma_2^2 x_2.
\end{array}  \right.  \nonumber  
\end{align}
\ $\square$
\end{corollary}
Figure \ref{fig:ex}.b shows an example of the cross-section of the Barrier surface for given positions of the pursuers. Similarly, the projection of the sets $\mathcal{R}_e$ and $\mathcal{R}_p$ into the 2D plane is also illustrated in Fig. \ref{fig:ex}.a; they are denoted, respectively, by $R_e$ and $R_p$.

\textit{Remark}. The same approach can be extended to consider more than two pursuers against \textit{one evader} in a straightforward manner. When more than one evader is considered, tractability calls for the assignment of pursuers to evaders to be performed in order to obtain the Barrier surface particular to that assignment. Assignment is necessary to prevent a pursuer trying to simultaneously pursue more than one evader which renders the herein constructed Barrier surface invalid and the solution of the attendant multi-player differential game intractable.


\section{Conclusions} \label{sec:concl}
The game of kind of the cooperative football differential game was solved: The Barrier surface that separates the teams' winning regions is constructed. Compared to previous work, the Barrier surface has been explicitly obtained in closed form. The cross-sections of the Barrier surface have been characterized in a simple way in terms of quadratic curves: circles and hyperbolae. The more general case of fast pursuers with different speeds has been addressed together with the classic case where all players, including the evader, have the same speed.

\bibliographystyle{IEEEtran}
\bibliography{ReferencesTAD}

\begin{thebibliography}{10}
\providecommand{\url}[1]{#1}
\csname url@samestyle\endcsname
\providecommand{\newblock}{\relax}
\providecommand{\bibinfo}[2]{#2}
\providecommand{\BIBentrySTDinterwordspacing}{\spaceskip=0pt\relax}
\providecommand{\BIBentryALTinterwordstretchfactor}{4}
\providecommand{\BIBentryALTinterwordspacing}{\spaceskip=\fontdimen2\font plus
\BIBentryALTinterwordstretchfactor\fontdimen3\font minus
  \fontdimen4\font\relax}
\providecommand{\BIBforeignlanguage}[2]{{%
\expandafter\ifx\csname l@#1\endcsname\relax
\typeout{** WARNING: IEEEtran.bst: No hyphenation pattern has been}%
\typeout{** loaded for the language `#1'. Using the pattern for}%
\typeout{** the default language instead.}%
\else
\language=\csname l@#1\endcsname
\fi
#2}}
\providecommand{\BIBdecl}{\relax}
\BIBdecl

\bibitem{Weintraub20}
I.~Weintraub, M.~Pachter, and E.~Garcia, ``An introduction to pursuit-evasion
  differential games,'' in \emph{arXiv preprint arXiv:2003:05013}, 2020.

\bibitem{Isaacs65}
R.~Isaacs, \emph{Differential Games}.\hskip 1em plus 0.5em minus 0.4em\relax
  New York: Wiley, 1965.

\bibitem{zhou2016cooperative}
Z.~Zhou, W.~Zhang, J.~Ding, H.~Huang, D.~M. Stipanovi{\'c}, and C.~J. Tomlin,
  ``Cooperative pursuit with voronoi partitions,'' \emph{Automatica}, vol.~72,
  pp. 64--72, 2016.

\bibitem{pierson2016intercepting}
A.~Pierson, Z.~Wang, and M.~Schwager, ``Intercepting rogue robots: An algorithm
  for capturing multiple evaders with multiple pursuers,'' \emph{IEEE Robotics
  and Automation Letters}, vol.~2, no.~2, pp. 530--537, 2016.

\bibitem{Selvakumar17}
J.~Selvakumar and E.~Bakolas, ``Evasion with terminal constraints from a group
  of pursuers using a matrix game formulation,'' in \emph{American Control
  Conference}, 2017, pp. 1604--1609.

\bibitem{Casini19}
M.~Casini, M.~Criscuoli, and A.~Garulli, ``A discrete-time pursuit-evasion game
  in convex polygonal environments,'' \emph{Systems and Control Letters}, vol.
  125, pp. 22--28, 2019.

\bibitem{chen2016multiplayer}
M.~Chen, Z.~Zhou, and C.~J. Tomlin, ``Multiplayer reach-avoid games via
  pairwise outcomes,'' \emph{IEEE Transactions on Automatic Control}, vol.~62,
  no.~3, pp. 1451--1457, 2017.

\bibitem{margellos2011hamilton}
K.~Margellos and J.~Lygeros, ``Hamilton--jacobi formulation for reach--avoid
  differential games,'' \emph{IEEE Transactions on Automatic Control}, vol.~56,
  no.~8, pp. 1849--1861, 2011.

\bibitem{Blake04}
M.~A. Blake, G.~A. Sorensen, J.~K. Archibald, and R.~W. Beard, ``Human assisted
  capture-the-flag in an urban environment,'' in \emph{IEEE International
  Conference on Robotics and Automation}, 2004, pp. 1167--1172.

\bibitem{Huang15}
H.~Huang, J.~Ding, W.~Zhang, and C.~J. Tomlin, ``Automation-assisted
  capture-the-flag: a differential game approach,'' \emph{IEEE Transactions on
  Control Systems Technology}, vol.~23, no.~3, pp. 1014--1028, 2015.

\bibitem{GarciaCDC18}
E.~Garcia, D.~W. Casbeer, and M.~Pachter, ``The capture-the-flag differential
  game,'' in \emph{57th IEEE Conference on Decision and Control}, 2018, pp.
  4167--4172.

\bibitem{fisac2015pursuit}
J.~F. Fisac and S.~S. Sastry, ``The pursuit-evasion-defense differential game
  in dynamic constrained environments,'' in \emph{IEEE 54th Annual Conference
  on Decision and Control}, 2015, pp. 4549--4556.

\bibitem{Coon17}
M.~Coon and D.~Panagou, ``Control strategies for multiplayer
  target-attacker-defender differential games with double integrator
  dynamics,'' in \emph{56th IEEE Conf. on Decision and Control}, 2017, pp.
  1496--1502.

\bibitem{Garcia2019}
E.~Garcia, D.~W. Casbeer, and M.~Pachter, ``Design and analysis of
  state-feedback optimal strategies for the differential game of active
  defense,'' \emph{IEEE Transactions on Automatic Control}, vol.~64, no.~2, pp.
  553--568, 2019.

\bibitem{harini2015new}
R.~H. Venkatesan and N.~K. Sinha, ``A new guidance law for the defense missile
  of nonmaneuverable aircraft,'' \emph{IEEE Transactions on Control Systems
  Technology}, vol.~23, no.~6, pp. 2424--2431, 2015.

\bibitem{Weintraub2018}
I.~E. Weintraub, E.~Garcia, and M.~Pachter, ``A kinematic rejoin method for
  active defense of non-maneuverable aircraft,'' in \emph{2018 American Control
  Conference}.\hskip 1em plus 0.5em minus 0.4em\relax IEEE, 2018, pp.
  6533--6538.

\bibitem{Shishika18}
D.~Shishika and V.~Kumar, ``Local-game decomposition for multiplayer
  perimeter-defense problem,'' in \emph{IEEE 57th Conference on Decision and
  Control}, 2018, pp. 2093--2100.

\bibitem{Garcia20LCS}
E.~Garcia, D.~W. Casbeer, and M.~Pachter, ``Optimal strategies of the
  differential game in a circular region,'' \emph{IEEE Control Systems
  Letters}, vol.~4, no.~2, pp. 492--497, 2020.

\bibitem{YanShi17}
R.~Yan, Z.~Shi, and Y.~Zhong, ``Defense game in a circular region,'' in
  \emph{IEEE 56th Conference on Decision and Control}.\hskip 1em plus 0.5em
  minus 0.4em\relax IEEE, 2017, pp. 5590--5595.

\bibitem{Oyler16}
D.~W. Oyler, P.~T. Kabamba, and A.~R. Girard, ``Pursuit--evasion games in the
  presence of obstacles,'' \emph{Automatica}, vol.~65, pp. 1--11, 2016.

\bibitem{Scott13}
W.~Scott and N.~E. Leonard, ``Pursuit, herding and evasion: A three-agent model
  of caribou predation,'' in \emph{American Control Conference}, 2013, pp.
  2978--2983.

\bibitem{GarciaACC19}
E.~Garcia, D.~W. Casbeer, A.~{Von Moll}, and M.~Pachter, ``Cooperative
  two-pursuer one-evader blocking differential game,'' in \emph{2019 American
  Control Conference}, 2019, pp. 2702--2709.

\bibitem{Bopardikar11}
S.~Bopardikar, S.~L. Smith, and F.~Bullo, ``On vehicle placement to intercept
  moving targets,'' \emph{Automatica}, vol.~47, no.~9, pp. 2067--2074, 2011.

\bibitem{Vonmoll19}
A.~{Von Moll}, E.~Garcia, D.~W. Casbeer, M.~Suresh, and S.~C. Swar, ``Multiple
  pursuer single evader border defense differential game,'' in \emph{AIAA
  Scitech 2019 Forum, AIAA 2019-1162}, 2019.

\bibitem{YanShi19}
R.~Yan, Z.~Shi, and Y.~Zhong, ``Task assignment for multiplayer reach-avoid
  games in convex domains via analytical barriers,'' \emph{IEEE Transactions on
  Robotics}, 2019.

\end{thebibliography}

\end{document}